\documentclass[11pt]{article}
\input{epsf.sty}
\usepackage{hyperref}
\usepackage{amsmath, amssymb}
\newtheorem {thm}{Theorem}[section]
\newtheorem {prop}[thm]{Proposition}
\newtheorem {lem}[thm]{Lemma}

\newtheorem {defn}[thm]{Definition}
\newtheorem {conj}[thm]{Conjecture}

\def\Cox{\hfill \Box}

\def\P{{\Bbb P}}
\def\Q{{\Bbb Q}}
\def\E{{\Bbb E}}

\def\0{{\bf 0}}

\def\b{\beta}

\def\d{\delta}

\def\e{\varepsilon}

\def\phi{\varphi}

\def\l{\lambda}

\def\o{\omega}

\def\L{\Lambda}

\def\T{\T}

\begin{document}

\title{On the Purity of the \\ 
free boundary condition Potts measure on random trees} 

\author{
Marco Formentin
\footnote{
Universit\`a degli Studi di Padova,
Dipartimento di Matematica Pura ed Applicata,
Via Trieste, 63 - 35121 Padova, 
Italia
\texttt{formen@math.unipd.it}}
\, and
Christof K\"ulske
\footnote{
University of Groningen, 
Department of Mathematics and Computing Sciences, 
Nijenborgh 9,   
9747 AC Groningen, 
The Netherlands
\texttt{kuelske@math.rug.nl},
\texttt{ http://www.math.rug.nl/$\sim$kuelske/ }}
}

\maketitle

\begin{abstract}
We consider the free boundary condition Gibbs measure 
of the Potts model on a random tree. 
We provide an explicit temperature interval below 
the ferromagnetic transition temperature for which this measure is extremal,   
improving older bounds of Mossel and Peres.
In information theoretic language extremality of the Gibbs measure
 corresponds to non-reconstructability for symmetric $q$-ary channels. 
The bounds for the corresponding threshold value of the inverse temperature 
are optimal for the Ising model and differ from the Kesten Stigum bound
by  only $1.50$ percent in the case $q=3$ and $3.65$ percent for $q=4$, independently of $d$. 
 Our proof uses an iteration of random boundary entropies from the outside of the tree to the inside, along with a symmetrization argument. 

\end{abstract}

\bigskip 


\section{Introduction} \label{sect:intro}

Interacting stochastic models on trees 
and lattices often differ in a fundamental way:
where a lattice model has a single transition point 
(a critical value for a parameter of the model), the corresponding 
model on a tree might possess multiple transition points. 
Such phenomena happen more generally for non-amenable graphs (where 
surface terms are no smaller than volume terms), trees 
being major examples \cite{Ly00}. 

A main example of an interacting model is the usual ferromagnetic Ising model. 
Here the interesting property which gives rise to a new transition 
is the extremality of the free b.c. (boundary condition) state. 
In an Ising model on the lattice, below the ferromagnetic transition temperature
the free boundary limiting measure will be a symmetric combination between 
the plus-state and the minus-state. On the tree, however, 
the open boundary state 
will still be extremal in a temperature interval strictly 
below the ferromagnetic transition temperature.  It ceases 
to be extremal at even lower temperatures. 

Ferromagnetic order on a tree is characterized by the fact 
that a plus-boundary condition at the leaves of a finite tree of depth $n$ 
persists to have influence on the origin when $n$ tends to infinity.
For the tree it now happens in a range of temperatures that, even though an all plus-boundary condition 
will be felt at the origin, a {\em typical boundary condition} chosen from
the free b.c. measure itself will not be felt at the origin for a range of temperatures below 
the ferromagnetic transition. The latter implies the extremality of the free b.c. state.   

We write throughout the paper 
$\theta = \tanh \b$ where $\b$ is the inverse temperature of the Ising (or Potts) model 
and denote by $d$ the number of children on a regular rooted tree. 
Then the Ising ferromagnetic transition temperature is given by   
$d \theta  =1$, and the transition temperature where the free b.c. state 
ceases to be extremal is given by $d \theta^2 =1$.  

A proof of the latter fact is contained in \cite{BleRuZa95}.
A beautiful alternate proof of the extremality for $d \theta^2 \leq 1$
for regular trees was given by Ioffe \cite{Io96a}. 
The method used therein was elegant but very much dependent 
on the two-valuedness of the Ising spin variable. This was exploited 
for the control of conditional probabilities in terms of projections to 
products of spins. Some care is necessary to treat the marginal case where 
equality holds in the condition. 
Indeed, one needs to control quadratic terms in a recursion; this is difficult for a general tree where 
the degrees are not fixed. A second paper 
\cite{Io96b} proves an analogue of the condition for general trees 
with arbitrary degrees but leaves this case open. 
Finally, for a general tree which  
does not possess any symmetries, \cite{PePe06} give a sharp criterion for extremality in terms of capacities.    
It remains an open problem to determine the extremal measures and the measure 
in the extreme decomposition of the open b.c. state for 
$d \theta^2 > 1$.

Let us remark that the problem of extremality of the open b.c. state 
is equivalent to the so-called reconstruction problem: 
We send a signal (a plus or a minus) 
from the origin to the boundary, 
making a prescribed error probability (that is related 
to the temperature of the Ising model) 
at every edge of the tree. In this 
way one obtains a Markov chain indexed by the tree.  
The reconstruction problem on a tree is said to be solvable, 
if the measure, obtained on the boundary at distance $n$ 
by sending an initial $+$, keeps a finite variational distance to the measure  
obtained by sending a $-$, as $n$ tends to infinity. Nonsolvability of reconstruction 
is equivalent to the extremality of the open b.c. state \cite{Mo01,JaMo04}. This is to  
say that there can be no transport of information along the tree 
between root and boundary, for typical signals. 

\subsection{The Potts model}

We denote by $T^N$ a finite tree rooted at $0$ of depth $N$. 
Then the {\em free b.c. Potts measure on $T^N$} is the probability distribution 
$\P^N$ that assigns to a configuration   
$\eta_{T^N}=(\eta(v))_{v\in T^N}\in \{1,2,\dots,q\}^{T^N}$  
the probability weights 
\begin{equation}\begin{split}\label{1}
&\P^N( \eta_{T^N} )= \frac{\exp\bigl( 
2\b \sum_{(v,w)}\d_{\eta(v),\eta(w)}
\bigr) }{Z_{\b,{T^N}}},
\cr
\end{split}
\end{equation}
where the sum is over all edges $(v,w)$ of the tree $T^N$ and $Z_{\b,{T^N}}$
is the partition function that makes the r.h.s. a probability measure. 

The {\em free b.c. Potts measure on an infinite tree $T$} 
is by definition the  ak limit  $\P= \lim_{N\uparrow \infty} \P_{T^N}$ 
when $T^N$ is an exhaustion of $T$. 
$\P$ is identical to what is called 
{\em the symmetric chain on $q$ symbols} in the context of the 
reconstruction problems in \cite{Mo01}. This chain 
has one parameter, namely the probability 
to change the symbol that is transmitted to any of the 
$q-1$ others, which is given by $\frac{1}{e^{2\b}+q-1}$. 

Recalling the DLR equations, {\em a Potts Gibbs measure} on a graph with vertex set 
$G$ is any measure $\P$ such that, for all finite subsets $\L\subset G$, the corresponding 
conditional probabilities of $\P$ are given by 
\begin{equation}\begin{split}\label{localspecification}
&\P( \eta_{\L}| \bar \eta_{G\backslash V} )= 
\frac{
\exp\bigl( 2\b \sum_{ {(v,w)}\atop{v,w \in \L} }
\d_{\eta(v),\eta(w)} + \sum_{{(v,w)}\atop {v\in \L, w \in G\backslash \L }}\d_{\eta(v),\bar\eta(w)}
\bigr) 
}
{
Z_{\b,{\L}}^{\bar \eta}
},
\cr
\end{split}
\end{equation}
where the sums are again along edges $(v,w)$ of the graph. 

Clearly the free b.c. measure on an infinite tree $T$ is a Gibbs measure. 
Recall that a Gibbs measure is said to be extremal if it can not be written as convex combination 
of other Gibbs measures. 

\subsection{Random trees}


Consider a random tree $T$ with vertices 
$i$ and number of children at the site $i$ given by $d_i$. We choose $d_i$ to be independent random variables 
with the same distribution $\Q$.  We use the symbol 
$\Q$ also  to describe the expected value. 
As is well known these appear as local approximations 
of random graphs which has newly emphasized 
their interest \cite{GeMo07}.
Our results however are already interesting in the 
case of regular trees where every vertex $i$ 
has precisely $d$ children.

\subsection{A criterion for extremality on random trees}

In this situation our main result, formulated for a random tree, is the following.  
Write $$P=\{(p_i)_{i=1,\dots,q}, p_i \geq 0\,\, \forall i,\,\, \sum_{i=1}^q p_i =1\}$$ for the simplex 
of Potts probability vectors. 

\begin{thm} \label{thm1} 
The free boundary condition Gibbs measure $\P$ 
is extremal, for $\Q$-a.e. tree $T$ when the condition    
$\Q(d_0) \frac{2\theta}{q-(q-2)\theta} \bar c(\b,q) <1$ is satisfied. 
Here, 
\begin{equation}\label{cbar}
\begin{split}
\bar c(\b,q):=
\sup_{p\in P }
\frac{\sum_{i=1}^q ( q p_i -1)\log (1+(e^{2 \b}-1) p_i)}{\sum_{i=1}^q ( q p_i -1)\log q p_i}\,.\cr
\end{split}
\end{equation}
\end{thm}

{\bf Remark:} It appears that the supremum over $P$ is achieved 
at the symmetric point $\frac{1}{q}(1,1,\dots,1)$ {\em only} in the Ising model $q=2$.  
This implies the sharpness of the bound in the Ising case, see also the discussion 
at the end of the paper.  
It is not surprising that the Potts model shows pecularities in comparison with the Ising model. 
That Potts is more intricate is seen already on the level of the much simpler problem 
of determining the ferromagnetic transition temperature (where the Gibbs measure 
becomes unique). Due to the lack of concavity of the r.h.s. of the recursion 
relation the transition is first (instead of second) order. 

{\bf Remark:} 
The best bound which has been previously given appears in \cite{MoPe03} on a $d$-ary tree.  
We recover it from our bounds when we use the estimate $\bar{c}(\beta,q)\leq \theta$ which will 
be discussed below.   
Moreover, numerically $\bar c(\beta,q)$ 
seems to decrease monotonically in $q$ at fixed $\beta$.  \bigskip\bigskip

Note also the bounds of Martinelli et al. \cite{Ma04} (see Theorem 9.3., Theorem 9.3.'
Theorem 9.3'') who give a nice criterion for non-reconstruction involving a Dobrushin constant 
of the corresponding Markov specification 
which however give worse estimates in the Potts model.

Let us put our result in perspective. For the purpose of the discussion 
we specialize to the case 
of the regular tree with $d$ children.
Denote by $\P^{N,k}$ the measures on $T^N$ obtained by putting 
the boundary condition $k$ to all Potts-spins at the outer boundary, and denote 
by $\P^k$ the corresponding limiting measures on $T$.  

Absence of ferromagnetic order (uniqueness of the Gibbs measure) 
can be detected by the fact that the distribution of the spin 
$\eta_0$ at the origin under the infinite volume measure $\P^k$ is the equidistribution, 
independently of the boundary condition $k$. 
This condition is easy to obtain by considering a simple 
one-dimensional recursion of numbers (instead of measures).  For more details 
see Subsection \ref{2.2}. 
Absence of ferromagnetic order in particular implies purity of the free b.c. state. 
In the language of the reconstruction problem this means non-solvability 
and as such the condition  is mentioned as Proposition 4 in  \cite{Mo01}.  

Let us compare with opposite results: 
It is known from the so-called Kesten-Stigum bound \cite{KeSt66}
that $d \l_2(\theta,q)^2> 1$ implies reconstructability 
(i.e. non-extremality of the free b.c. measure). 
Here $\l_2(\theta,q)$ is the second eigenvalue of the transition matrix 
that produces the free b.c. Potts model by broadcasting from the origin to the boundary;  
it is decreasing in $q$ at fixed $\theta$, 
and  increasing in $\theta$ at fixed $q$. This is intuitively clear: the bigger the number of states $q$ and 
the smaller the inverse temperature, the easier  it is 
to forget about the information put at the boundary.
Moreover it is proved as Theorem 2 in \cite{Mo01}
that when one fixes $d$ and a value of $d \l_2(\theta,q)\equiv \l>1$,  
for $q$ large enough the reconstruction problem is solvable for the corresponding 
value of $\theta$. 

Now, our method of proof is based on controlling 
recursions for the probability distributions at roots of subtrees 
from the outside to the inside of a tree. 
These are recursions on log-likelihood ratios of Potts 
probability vectors for the root of subtrees, and these ratios  
are random w.r.t. the boundary condition (which is chosen according to 
the free b.c. condition measure). 

Understanding recursions for probability distributions (needed to investigate 
the purity of the free b.c. state) is 
much less straightforward than controlling recursions for real numbers (needed 
for investigating the existence of ferromagnetic order). 
We prove convergence to a Dirac-distribution by controlling the boundary relative entropy, 
generalizing from the approach of \cite{PePe06} for the Ising model. 
Novelties appear for the Potts model, a key point being proper 
symmetrization 
to bring out the constant \eqref{cbar},  beginning with Lemma \ref{lemma2.99}.


\bigskip 

\textbf{Acknowledgements: }\\
The authors thank Aernout van Enter for interesting discussions and useful comments 
on the manuscript. 

\section{Proof}

To show the triviality of a measure $\mu$ 
on the tail sigma-algebra  it suffices to show that, 
for any fixed cylinder event $A$ we have 
\begin{equation}\label{2.1}
\begin{split}
\lim_{N\uparrow \infty}
\mu \left | \mu(A| {\cal T}_{N}) -\mu(A) \right |   =0,
\end{split}
\end{equation}
where ${\cal T}_{N}$ is the sigma-algebra 
created by the spins that have at least distance 
$N$ to the origin (see \cite{Ge88} Proposition 7.9).

We denote by $T^N$ the tree rooted at $0$ of depth $N$. 
The notation $T^N_v$ indicates the sub-tree of $T^N$ rooted at $v$ obtained from ``looking to the outside" on the tree $T^N$.  
We denote by $\mathbb{P}^{N,\xi}_v$ the correponding Potts-Gibbs measure on 
$T^N_v$ with boundary condition on 
$\partial T^N_v$ given by $\xi=(\xi_i)_{i\in \partial T^N_v}$. 
We denote by $\mathbb{P}^{N}_v$ the correponding Potts-Gibbs measure on 
$T^N_v$ with free boundary conditions, as in \eqref{1}.

We are going to show that the distribution of the probabilities to see a value $s$ 
at the origin, obtained by putting a boundary condition $\xi$ at distance $N$ that 
is chosen according to the free measure $\P$ itself, converges to the equidistribution
in probability. This reads  
\begin{equation}\label{2.2666}
\begin{split}
\lim_{N\uparrow \infty}
\P\Bigl(\xi :  \Bigl |\mathbb{P}^{N,\xi}(\eta(0)=s)- \frac{1}{q} \Bigr | \geq \e\Bigr) \rightarrow 0 .
\end{split}
\end{equation}
This then implies (\ref{2.1}). 


To achieve (\ref{2.2666}) it is more convenient 
to look at the probability distribution for the spin at the root $v$ obtained with the boundary 
condition $\xi$ in terms of the ``log-likelihood ratios" defined by 
\begin{equation}
\begin{split}
X^j_k(v;\xi):=\log\frac{\mathbb{P}^{N,\xi}_v(\eta(v)=j)}{\mathbb{P}^{N,\xi}_v(\eta(v)=k)} ,
\end{split}
\end{equation}
where $1\leq j \neq k  \leq q$. 
Ultimately we are interested to show the convergence of 
these quantities at $v=0$ to zero, for all pairs $j,k$, in $\P$-probability,   
as the depth $N$ of the tree tends to infinity.

We denote the measure at the boundary at distance $N$ from the root on the tree 
emerging from $v$, which is obtained by 
conditioning the spin in the site $v$ to take the value to be $j$, by 
\begin{equation}
\begin{split}
Q^{N,j}_v(\xi):=\mathbb{P}^N_v(\eta:\eta_{|\partial T_v^N}=\xi|\:\eta(v)=j) .
\end{split}
\end{equation}

\begin{defn} Denote the relative entropy of the boundary measures between the states 
obtained by conditioning the spin at $v$ to be $1$ respectively $2$, by  
\begin{equation}
\begin{split}
m_v^{(N)}=S(Q_v^{N,2}| Q_v^{N,1})=\int Q_v^{N,2}(d\xi)\log \frac{Q_v^{N,2}(\xi)}{Q_v^{N,1}(\xi)}.
\end{split}
\end{equation}
\end{defn}

Here and in the sequel denote by $w$ the children of $v$, indicated by the symbol 
$v\rightarrow w$.

\begin{lem} \label{lemma2.99} The boundary relative entropy can be written 
as an expected value w.r.t. the open boundary condition Gibbs measure $\P$ in the form 
\begin{equation}
\begin{split}
S(Q_v^{N,2}| Q_v^{N,1})=\frac{1}{q-1}\int \P(d\xi)\sum_{i=1}^q \phi\left( 
q \mathbb{P}^{N,\xi}_v(\eta(v)=i)
\right) ,
\end{split}
\end{equation}
with $\phi(x)=(x-1)\log x$. 
\end{lem}
\bigskip 

{\bf Proof: } 
In the first step we express the relative entropy as an expected value 
\begin{equation}
\begin{split}
S(Q_v^{N,2}| Q_v^{N,1})=q\int \P(d\xi)g\Bigl( \mathbb{P}^{N,\xi}_v(\eta(v)=2),
\mathbb{P}^{N,\xi}_v(\eta(v)=1)\Bigl),
\end{split}
\end{equation}
with 
\begin{equation}
\begin{split}
g(p_2,p_1)=p_2\log\frac{p_2}{p_1}.
\end{split}
\end{equation}
To see this, we use that 
\begin{equation}
\begin{split}
\frac{dQ^{N,2}_v}{d\mathbb{P}_v^N}(\xi)=q \mathbb{P}^{N,\xi}_v(\eta(v)=2), 
\end{split}
\end{equation}
by the definition of the conditional probability and the fact that the marginal 
of $\P$ at any site is the equidistribution. 

In the next step we use the invariance of $\P$ under permutation of the Potts-indices  to write 

\begin{equation}
\begin{split}
S(Q_v^{N,2}| Q_v^{N,1})=q\int \P(d\xi)(R g)\Bigl(\mathbb{P}^{N,\xi}_v(\eta(v)=1), \mathbb{P}^{N,\xi}_v(\eta(v)=2),\dots,
\mathbb{P}^{N,\xi}_v(\eta(v)=q)\Bigl),
\end{split}
\end{equation}
where $R$ is the symmetrization operator acting on functions $f(p_1,\dots,p_q)$ of Potts-probability vectors by 
\begin{equation}\label{Roperator}
\begin{split}
(R f)(p_1,p_1,\dots,p_q)=\frac{1}{q!}\sum_{\pi} f(p_{\pi(1)}, p_{\pi(2)},\dots,p_{\pi(q)}),
\end{split}
\end{equation}
where $\pi$ runs over the permutations of $\{1,\dots,q\}$.

One verifies that 
\begin{equation}
\begin{split}
(R g)(p_1,p_1,\dots,p_q)=\frac{1}{q(q-1)}\sum_{i=1}^q (q p_i -1)\log q p_i ,
\end{split}
\end{equation}
which proves the lemma.$\Cox$

\subsection{Recursions for the boundary entropy for subtrees}

\begin{prop}\label{prop1} The boundary relative entropy $m_v^{(N)}$ at the site $v$ obeys 
the following linear 
recursive inequalities in terms of the values at the children $w$, given by 
\begin{equation}
\begin{split}
m_v^{(N)}\leq \frac{2\theta}{q-\theta(q-2)}\bar c(\b,q)\sum_{w:v\rightarrow w} m_w^{(N)}.
\end{split}
\end{equation}
\end{prop}

{\bf Remark:} Noting that $\frac{Q_v^{N,j}(\xi)}{Q_v^{N,k}(\xi)}=X^j_k(v;\xi)$ we may write 
\begin{equation}
\begin{split}
m_v^{(N)}=\int Q_v^{N,2}(d\xi) X_{1}^2 (v; \xi).
\end{split}
\end{equation}

{\bf Remark:} Suppose that we are considering a spherically symmetric 
tree. This means that the number of offspring depends only on the generation, e.g. 
$d_v=d_{|v|}$ where $|v|$ is the distance of $v$ to the origin (that is 
the length of the unique path from the origin to $v$). 
Then $m^{(N)}_{v}=m^{(N)}_{|v|}$ and so 
\begin{equation}
\begin{split}
m_{k}^{(N)}\leq \frac{2\theta}{q-\theta(q-2)}\bar c(\b,q)d_k m_{k+1}^{(N)}.
\end{split}
\end{equation}
So $\lim_{N\uparrow \infty}m_{0}^{(N)}=0$ is implied by 
$\sum_{k=1}^\infty \log (c d_k)=-\infty$ with $c= \frac{2\theta}{q-\theta(q-2)}\bar c(\b,q)$. 

\bigskip 
{\bf Proof of Theorem 1.1}
Taking expectation w.r.t. the random tree
we note that $\E m_v^{(N)}=\E m_{|v|}^{(N)}$. Now, 
using Wald's inequality we have 
\begin{equation}
\begin{split}
\E m_k^{(N)}\leq \frac{2\theta}{q-\theta(q-2)}\bar c(\b,q)
\E d_0 \,\E(m_{k+1}^{(N)}) .
\end{split}
\end{equation}
From this follows that $\lim_{N\uparrow\infty}\E m_{0}^{(N)}=0$
using the uniform boundedness in $N$,  
$\E m_{N-1}^{(N)}\leq C\E(d_0)$. This can be seen from   
Lemma \ref{lemma2.3} a few lines below.   
 $\Cox$

\bigskip

To prove Proposition \ref{prop1} at first a recursion for the log-likelihood ratios 
$X^j_k(v;\xi)$ has to be derived, for fixed finite tree of depth $N$ from 
the outside to the inside. This iteration is standard, but we include its 
derivation for the convenience of the reader. 
In the following we omit the dependence on the fixed boundary 
condition $\xi$ in the notation.

\begin{lem}\label{lemma2.3} For all indices $1\leq j, k  \leq q$ we have 
\begin{equation}
\begin{split}\label{iteration}
X^j_k(v)=\sum_{\o:v\rightarrow w}\log\frac{\sum_{i\neq k,j}\exp[X_k^i(w)]+1+\exp(2\beta)
\exp[X_k^j(w)]}{\sum_{i\neq k,j}\exp[X_k^i(w)]+\exp(2\beta)+\exp[X_k^j(w)]} .
\end{split}
\end{equation}
\end{lem}

{\bf Proof: }
 Note that the Potts-measure 
 $\mathbb{P}^{N,\xi}_v$ is proportional to the weight 
 $$W(\eta)=\prod_{x\rightarrow y, x\geq v}\exp[2\beta \delta_{\eta(x),\eta(y)}],$$ 
 where the product is taken over the neighboring vertices coming after $v$ looking from the root of the tree. The normalization factor will be $Z_v^{-1}$.\\
We want to rewrite $X^j_k(v)$ as a function of $X^j_k(w)$ where $w$ are the children of $v$.
The key observation is that 
$$W(\eta_v)=\prod_{w :v\rightarrow w}W(\eta_{w})\exp[2\beta \delta_{\eta(v),\eta(w)}],$$ where we have written $\eta_{v}$ for the restriction of $\eta$ to the sub-tree $T_{v}^N$. Now,
\begin{equation}
\begin{split}
&\mathbb{P}^{N,\xi}_v(\eta(v)=j)=Z_v^{-1}\prod_{w:v\rightarrow w}\sum_{\eta_{w}}W(\eta_{w})\exp[2\beta \delta_{j,\eta(w)}]\cr
&=Z_v^{-1}\prod_{w:v\rightarrow w}Z_{w}\sum_{i=1}^q Z_{w}^{-1}\exp[2\beta \delta_{j,i}]\sum_{\eta_{w}:\eta(w)=i}W(\eta_{w})\cr
&=Z_v^{-1}\prod_{w:v\rightarrow w}Z_{w}\sum_{i=1}^q\exp[2\beta \delta_{j,i}]\,\mathbb{P}^{N,\xi}_{w}(\eta(w)=i).
\end{split}
\end{equation}
The same computation can be done for $\mathbb{P}^{N,\xi}_v(\eta(v)=k)$ to obtain:
$$\mathbb{P}^{N,\xi}_v(\eta(v)=k)=Z_v^{-1}\prod_{w:v\rightarrow w}Z_{w}\sum_{i=1}^q\exp[2\beta \delta_{k,i}]\mathbb{P}^{N,\xi}_{w}(\eta(w)=i).$$ 
Now consider the ratio and then divide everything by $\mathbb{P}^{N,\xi}_{w}(\eta(w)=k)$:
$$\frac{\mathbb{P}^{N,\xi}_v(\eta(v)=j)}{\mathbb{P}^{N,\xi}_v(\eta(v)=k)}=
\prod_{w:v\rightarrow w}\frac{\sum_{i=1}^q\exp[2\beta \delta_{j,i}]\mathbb{P}^{N,\xi}_{w}(\eta(w)=i)}{\sum_{i=1}^q\exp[2\beta \delta_{k,i}]\mathbb{P}^{N,\xi}_{w}(\eta(w)=i)}=$$
$$=\prod_{w:v\rightarrow w}\frac{\sum_{i\neq k,j}\frac{\mathbb{P}^{N,\xi}_{w}(\eta(w)=i)}{\mathbb{P}^{N,\xi}_{w}(\eta(w)=k)}+1+
\exp(2\beta)\frac{\mathbb{P}^{N,\xi}_{w}(\eta(w)=j)}{\mathbb{P}^{N,\xi}_{w}(\eta(w)=k)}
}{\sum_{i\neq k,j}\frac{\mathbb{P}^{N,\xi}_{w}(\eta(w)=i)}{\mathbb{P}^{N,\xi}_{w}(\eta(w)=k)}+\exp(2\beta)+
\frac{\mathbb{P}^{N,\xi}_{w}(\eta(w)=j)}{\mathbb{P}^{N,\xi}_{w}(\eta(w)=k)}},$$
which proves the result.$\Cox$

\subsection{The ferromagnetic ordering}\label{2.2}

Let us quickly deviate from the proof of Proposition \ref{prop1} and 
discuss the threshold value for the ferromagnetic ordering (where the infinite 
volume states with uniform boundary conditions cease to be different). 

Observe that for a boundary condition $\xi$ that is all $q$ we have 
that $X^j_k(v)=0$ for all $1\leq i,j\leq q-1$, and further  that 
$X_i^q(v) =X_1^q(v)$ for all $i=1,\dots,q-1$. 
So the iteration runs on the one-dimensional quantity $X_1^q(v)$ and reads 
\begin{equation}
\begin{split}
X^q_1(v)&=\sum_{\o:v\rightarrow w}\log\frac{q-1+\exp(2\beta)
\exp[X_1^q(w)]}{q-2+\exp(2\beta)+\exp[X_1^q(w)]} \cr
&=:\sum_{\o:v\rightarrow w}\psi(X_1^q(w)).
\end{split}
\end{equation}
For a regular tree with $d$ children we have 

\begin{equation}
\begin{split}
X^q_1(k)&=d \psi(X_1^q(k+1)).
\end{split}
\end{equation}

We have to distinguish now the cases of $q=2$ and $q\geq 3$. 
For $q=2$ we see by computation of the second derivative 
that the function $\psi$ is concave. This means that 
the critical value $\b$ for which a positive solution $X$ ceases 
to exist is given by $1=d \psi'(0)$. 

The derivative at $X=0$ (which we state now for general $q$) reads 
\begin{equation}
\begin{split}
\frac{\partial}{\partial X}\psi(X)\bigl |_{X=0}= \frac{e^{2 \b}-1}{e^{2 \b}+q -1}=\frac{2 \theta}{q - (q-2)\theta}. 
\end{split}
\end{equation}
Hence, the critical value in the Ising case is given by $d \tanh \b=1$, 
for a regular tree where every vertex has $d$ children. 

We note that this quantity equals $\l_2$, the second eigenvalue of the 
transition matrix associated to the model. 

Let us now turn to the Potts model with $q\geq 3$. 
A computation shows that $\psi''(0)>0$ for $\b>0$ and $q\geq 3$, 
and hence the function $\psi$ is {\em not } concave. This reflects 
the fact that the transition at the critical point where a positive solution 
ceases is a first order transition, where the nonzero solution is bounded 
away from zero. 

For a regular tree with $d$ children we can derive the 
transition value $\b(q,d)$ as follows: 
We must have $1=d\psi'(X^*)$, meaning that the function $\psi$ 
touches the line $X$ with the same slope.  This equation translates into 
$\frac{1}{d}=\frac{a x}{q-1+a x}- \frac{x}{q-2+a+ x}$ in the variables $a=e^{2 \b}$, $x=\exp[X^*]$. 
The fixed point equation itself reads $x^\frac{1}{d}= 
\frac{q-1+a x}{q-2+a+x}$ . 

From these two equations the critical values can be derived numerically 
for any $d,q$. We note moreover that, for the special case of a binary tree $d=2$,  the 
fixed point equation is cubic in the variable $y:=x^\frac{1}{2}$.
 The fixed point equation is equivalent to $y (q - 2 + a + y^2) - ((q - 1) + a y^2)=0$. 
We already know one root, it is $y=1$, so we can produce a quadratic equation 
by polynomial division. Writing $y=1+u$ we get the solutions  
$u=\frac{1}{2}(-3 + a - \sqrt{5 - 2 a + a^2 - 4 q})$ and 
$u=\frac{1}{2}(-3 + a + \sqrt{5 - 2 a + a^2 - 4 q})$. The solution ceases to exist 
when the argument of the squareroot becomes negative which 
results in a critical value $a=1 + 2 \sqrt{q-1}$, or $\b(d=2,q)=\frac{1}{2}\log (1 + 2 \sqrt{q-1})$.
We note the numerical values $\b(d=2,q=3)=0.671227$, $\b(d=2,q=4)=0.748034$.

The same type of reasoning can be used for $d=3$ where the fixed point equation 
requires the solution of a fourth order equation in $z=x^\frac{1}{3}$, 
which can be reduced to a third order equation by dividing out the root $z=1$.
We don't give details here.

\subsection{Controlling the recursion relation for the boundary entropy}

\begin{lem}\label{lemma2.5}
\begin{equation}
X^j_i(v)=\sum_{\o:v\rightarrow w}
\Bigl[u\Bigl(\mathbb{P}^{N,\xi}_v(\eta(v)=j)\Bigr) 
- u\Bigl(\mathbb{P}^{N,\xi}_v(\eta(v)=i)
\Bigr)\Bigr] ,
\end{equation}
where 
\begin{equation}
u( p_1) =\log(1+p_1(e^{2 \b}-1)) .
\end{equation}
\end{lem}

{\bf Proof:  }
Remember the recursion given in Lemma \ref{lemma2.3}. 
Now re-express the $X$'s by the $p$-variables and use the fact that they 
form a probability vector.$\Cox$ 
\bigskip

Using this we may derive the following equality
on the iteration of the boundary entropy. 

\begin{lem}
\begin{equation}\label{meq}
\begin{split}
Q^{N,2}_v X^2_1(v)&=\frac{2\theta}{q-(q-2)\theta}\sum_{\o:v\rightarrow w}
Q^{N,2}_{w }
\Bigl[u\Bigl(\mathbb{P}^{N,\xi}_w(\eta(w)=2)\Bigr) 
- u\Bigl(\mathbb{P}^{N,\xi}_w(\eta(w)=1)
\Bigr)\Bigr]  .  \cr 
\end{split}
\end{equation}
\end{lem}

{\bf Proof}: As the second piece of information next to Lemma \ref{lemma2.5} which is needed to understand 
the iteration for the boundary relative entropy $m_v^{(N)}$ 
we must see how the boundary measure $Q^{N,j}_v(d\xi)$, obtained by conditioning 
at $v$, relates to the boundary measures obtained by conditioning 
at the  children, denoted by $w$. 

For the Potts model a computation shows that 
\begin{equation}\label{eq21}
\begin{split}
Q^{N,j}_v&=\prod_{v\rightarrow w}\left[\frac{\exp(2\beta)}{(q-1)+\exp(2\beta)}Q^{N,j}_{w}+\frac{1}{(q-1)+\exp(2\beta)}\sum_{i\neq j}Q^{N,i}_{w}\right]\cr 
&=\prod_{v\rightarrow w}\left[\frac{1+\theta}{q-(q-2)\theta}
Q^{N,j}_{w }+\frac{1-\theta}{q-(q-2)\theta}\sum_{i\neq j}Q^{N,i}_{w}\right].\cr 
\end{split}
\end{equation}

Thus, to control the iteration we must look at the terms 

\begin{equation}\label{eq21}
\begin{split}
&\left[\frac{1+\theta}{q-(q-2)\theta}
Q^{N,2}_{w }+\frac{1-\theta}{q-(q-2)\theta}Q^{N,1}_{w}+
\frac{1-\theta}{q-(q-2)\theta}\sum_{i\geq 3}Q^{N,i}_{w}\right]\cr
&\Bigl[u\Bigl(\mathbb{P}^{N,\xi}_w(\eta(w)=2)\Bigr) 
- u\Bigl(\mathbb{P}^{N,\xi}_w(\eta(w)=1)
\Bigr)\Bigr]. \cr 
\end{split}
\end{equation}
We first note that, by symmetry under 
the measure  $Q^{N,i}_{w}$, for 
$i=3,\dots,q$, the corresponding terms in the sum vanish. 
Now we use the permutation symmetry of the Potts indices to see the proof.
$\Cox$

\bigskip 
Next we use the following representation.  

\begin{lem}
\begin{equation}\label{meq}
\begin{split}
Q^{N,2}_v X^2_1(v)&=\frac{2\theta}{q-(q-2)\theta}\sum_{w:v\rightarrow w}
\int \P(d\xi)h\Bigl( \mathbb{P}^{N,\xi}_w(\eta(w)=2),
\mathbb{P}^{N,\xi}_w(\eta(w)=1)\Bigl),
\end{split}
\end{equation}
\end{lem}
with 
\begin{equation}
\begin{split}
h(p_2,p_1)= q p_2 (u(p_2)-u(p_1)).
\end{split}
\end{equation}

{\bf Proof: } 
This follows as in the Proof of Lemma \ref{lemma2.99} by plugging in the Radon-Nikodym derivative  
of $Q^{N,2}_{w }$ w.r.t. the open b.c. measure. 
$\Cox$

\bigskip 

With these preparations we can now finish the proof of the main proposition. \bigskip 

{\bf Proof of Proposition \ref{prop1}: } Recalling the definition of the symmetrization operator \eqref{Roperator} 
we obtain 
\begin{equation}\label{meq}
\begin{split}
Q^{N,2}_v X^2_1(v)&=\frac{2\theta}{q-(q-2)\theta}\sum_{w:v\rightarrow w}
\int \P(d\xi)(R h)\Bigl( \mathbb{P}^{N,\xi}_w(\eta(w)=1),\dots, \mathbb{P}^{N,\xi}_w(\eta(w)=q)
\Bigr),
\end{split}
\end{equation}
where 
\begin{equation}\label{meq}
\begin{split}
(R h)(p_1,\dots,p_q)=\frac{1}{q-1}\sum_{i=1}^q ( q p_i -1)u( p_i).
\end{split}
\end{equation}
From here follows that 
\begin{equation}\label{HH}
\begin{split}
Q^{N,2}_v X^2_1(v)&=\frac{2\theta}{q-(q-2)\theta}\sum_{w:v\rightarrow w}
\int \P(d\xi)H\Bigl( \mathbb{P}^{N,\xi}_w(\eta(w)=1),\dots, \mathbb{P}^{N,\xi}_w(\eta(w)=q)
\Bigr),
\end{split}
\end{equation}
where 
\begin{equation}\label{meq}
\begin{split}
H(p_1,\dots,p_q)=\frac{1}{q-1}\sum_{i=1}^q ( q p_i -1)\tilde u( p_i),
\end{split}
\end{equation}
with
\begin{equation}
\tilde u( p_1) =\log\frac{1+p_1(e^{2 \b}-1)}{1+\frac{1}{q}(e^{2 \b}-1)} .
\end{equation}

From \eqref{HH}  we have the linear recursion relation 
\begin{equation}\label{HH1}
\begin{split}
&m^{N}(v)=Q^{N,2}_v X^2_1(v)\cr
&\leq\frac{2\theta}{q-(q-2)\theta} \bar c(\b,q) \sum_{w:v\rightarrow w}
\int \P(d\xi)R g\Bigl( \mathbb{P}^{N,\xi}_w(\eta(w)=1),\dots, \mathbb{P}^{N,\xi}_w(\eta(w)=q)
\Bigr)\cr
&\leq\frac{2\theta}{q-(q-2)\theta} \bar c(\b,q) \sum_{w:v\rightarrow w} m^{N}(w)
\end{split}
\end{equation}
and from here the result of the proposition follows.$\Cox$ 
\bigskip 
\bigskip

We could end the paper at this point, but let us comment 
on the constant appearing, and provide the following conjecture. 

Define
\begin{equation}
\begin{split}
\hat c(\b,q)
:=\sup_{p\in P, p_2=\dots =p_q}\frac{H(p_1,\dots,p_q)}{R g(p_1,\dots,p_q)}.
\end{split}
\end{equation}

\begin{conj}\label{lem2.10}
We believe that $\hat c(\b,q)=\bar c(\b,q)$. 
\end{conj}

We checked this numerically for small values of $q$. 
If the previous conjecture is true, the two properties of $\hat c(\b,q)$, 
namely, monotonicity in $q$ and the bound $\hat c(\b,q)\leq \theta$ carry over. 
These two properties are seen as follows. 

\begin{lem}\label{lem2.11}
\begin{equation}\label{cbar}
\begin{split}
\hat c(\b,q)=\sup_{x\in D_q} \bar{\varphi}(q,\lambda_q)(x),\cr
\end{split}
\end{equation}
with the function 
\begin{equation}
\bar{\varphi}(q,\lambda_q)(x)=\frac{\log\left(\frac{1+\lambda_q x}{1-\lambda_q (q-1)x}\right)}{\log\left(\frac{1+q x}{1-q(q-1)x}\right)},
\end{equation}
with parameter $\l_q=\frac{e^{2 \b}-1}{1+ \frac{1}{q}(e^{2 \b}-1)}$ on the range 
$D_q=\left[-\frac{1}{q},\frac{1}{q(q-1)}\right]$ with $D_{(q-1)}\supset D_q$.
\end{lem}

{\bf Proof: }  
Change to new coordinates on the simplex of probability vectors $(p_1,\dots,p_q)$ given by 
\begin{equation}\label{HH2}
\begin{split}
x_i&=p_i - \frac{1}{q} \text{ for }  i=1,\dots, q-1, \cr
\end{split}
\end{equation}
take $x=x_i$ for $i=1,\dots,q-1$, 
and use Lemma \ref{lem2.10}. 
$\Cox$

\begin{lem}
For all $q\geq 3$ we have that 
\begin{equation}
\hat{c}(\beta,q)<\hat{c}(\beta,q-1)\leq \frac{\lambda_2}{2}=\theta.
\end{equation}
\end{lem}

{\bf Proof:} 
We use that 
\begin{equation}
\frac{\partial \bar{\varphi}(q,\lambda_q)(x)}{\partial q}<0,
\end{equation}
for $x\in D_q$. This gives 
\begin{equation}\label{inequality}
\bar{\varphi}(q,\lambda_q)(x)<\bar{\varphi}(q-1,\lambda_{q-1})(x)<\bar{\varphi}(2,\lambda_2)(x),\:\:\:\:\:\: x\in D_q.
\end{equation}
$\Cox$

\noindent This observation makes it very simple to compute 
$\hat{c}(\beta,q)$ numerically, for every $q$.
\bigskip

Next, what about the sharpness of the constant? 
Could  it be possible that Theorem \ref{thm1} in fact holds with the sharp value 
$\frac{e^{2\b}-1}{q-1+ e^{2\b}}$ replacing the constant $\bar c(\b,q)$? 
In our approach such a conjecture would be based on looking at the Hessian of the function 
$$\partial_{x_i,x_j}\phi_{\l,q}(x_1,\dots,x_{q-1})\bigl |_{x_k=0\,\, \forall k}= 
4 \l  1_{i=j} + 2 \l  1_{i\neq j}.$$ Indeed, heuristically it should suffice 
to look at the quadratic approximation around the equidistribution. This results 
in the rigorous lower bound    
$\bar c(\b,q)\geq \frac{\l}{q}=\frac{e^{2\b}-1}{q-1+ e^{2\b}}$ which 
we recognize as the Kesten-Stigum bound. 
For the Ising model 
we have equality, which is not true for $q=3$.  

Let us compare with the recent literature. In their paper \cite{MeMo06} 
Montanari and Mezard make the conjecture 
that the Kesten-Stigum bound is sharp for  $q\leq 4$, or more precisely: 

\begin{conj}\label{1} (M\'ezard and Montanari 2006)
Consider the Potts model with $q$ symbols on a $d$-ary tree 
and let $\l_2=\frac{e^{2\b}-1}{e^{2\b}+q-1}=\frac{2\theta}{q-(q-2)\theta}$, with  $\theta=\tanh(\b)$, then, 
if $q\leq 4$ and $d<d_{max}$, there is reconstruction if and only if $d\l^2_2>1$.  
\end{conj}
This conjecture is based on extensive numerical simulations of the random recursion. 
Moreover, the restriction on $d$ comes from the limitation on the values of $d$ they can treat numerically
and they actually think that $d_{max}=+\infty$.
How close are the Kesten-Stigum bounds and our constants? 
We obtain  numerically 
$\bar c(\b,q)= \frac{e^{2\b}-1}{q-1+ e^{2\b}} (1+\e(q)) $ with $\e(3)=0.0150$ and $\e(4)=0.0365$.   
If we specialize to a binary tree, and take 
advantage of the possible temperature dependence of $\e$ we obtain 
$\b_c:=\sup\{\b :        2\frac{2\theta}{3-\theta} \bar c(\b,3)              <1 \} =1.0434$ for $q=3$ 
and $\b_c:=\sup\{\b :       2\frac{2\theta}{4-2\theta} \bar c(\b,4)              <1 \} =1.1555 $ in the case $q=4$.

After completion of the first draft of our present work Sly's preprint  \cite{Sly08} appeared where he proves the following. 
\begin{thm} (Sly 08)
When $q\leq 3$, and $d>d_{min}$, then the Kesten-Stigum bound is sharp,
while the Kesten-Stigum bound is never sharp when $q\geq 5$.
\end{thm} 

His method uses large degrees to justify 
quadratic expansions by means of central limit theorem approximation 
and makes no statements for small degrees where our estimates also apply. 
In view of these results it is natural for us to conjecture the following. 

\begin{conj} For $q\leq 4$ there is $\mathbb{Q}$-a.s. no reconstruction if  
\begin{equation}
\mathbb{Q}(d_0)\left(\frac{2\theta}{q-(q-2)\theta}\right)^2<1
\end{equation}
\end{conj}

Finally, what can we say about $q\geq 5$? 
Montanari and M\'ezard \cite{MeMo06}
find in all the test-cases of $q\geq 5$ and $d$ which they 
treat by simulations that the Kesten-Stigum bound is {\em not} sharp. Let us therefore conclude the paper by making 
a comparison of our and their values. Table 1 contains the simulation values from 
Montanari and M\'ezard for the critical error theshold $\e_r$ (probability to switch to a new symbol), taken from 
\cite{MeMo06}, 
as well as the corresponding numerical values of the critical inverse temperature $\b_r$ and the second eigenvalue $\l_r$. 
(The three columns are equivalent but we give them for easy comparison.)  
Table 2 contains as $\b_c$ our lower bound on the presumed true inverse 
reconstruction temperature $\b_r$, and the corresponding numerical value $\l_c$ of the second eigenvalue. 
We remark that our bounds appear to be close also here.

   \begin{table}
   \begin{tabular}{|c|c|c|c|} 
   \hline 
   q=5&$\epsilon_r$&$\beta_r=-0.5\log\left(\frac{\epsilon_r}{(q-1)(1-\epsilon_r)}\right)$&$\lambda_r=1-\frac{q}{q-1}\epsilon_r$\\
   \hline
   $d=2$&0.2348&1.2838&0.7065\\
   \hline
   $d=3$&0.33881&1.0285&0.5765\\
   \hline 
   $d=4$&0.4008&0.8942&0.4990 \\
   \hline
   $d=7$&0.4986&0.6955&0.3767 \\
   \hline
   $d=15$&0.5955&0.4998&0.2556 \\
   \hline
   \end{tabular}
   \bigskip\\
   \caption{Simulation results for the exact reconstruction thesholds  by Mezard and Montanari \cite{MeMo06}}
   \label{}
   \end{table}\bigskip 
   
   \begin{table}
   \begin{tabular}{|c|c|c|} 
   \hline 
   q=5&$\beta_c$&$\lambda_c$\\
   \hline
     $d=2$&1.2425&0.6875\\
   \hline
   $d=3$&0.98535&0.5526\\
   \hline 
   $d=4$&0.8520&0.47346 \\
   \hline
   $d=7$&0.65465&0.35095\\
   \hline
   $d=15$&0.4640&0.2342\\
   \hline
   \end{tabular}
   \bigskip\\
   \caption{Our bounds deduced from Theorem \ref{thm1} }
   \label{}
   \end{table}


\begin{thebibliography}{}

\bibitem{BleRuZa95} P. M. Bleher, J. Ruiz, V.A. Zagrebnov, 
{\em On the purity of limiting Gibbs state for the Ising model on the Bethe lattice, } 
J. Stat. Phys. \textbf{79}, 473-482 (1995).  


\bibitem{Ge88} H. O. Georgii, {\em 
Gibbs measures and phase transitions,} de Gruyter Studies in Mathematics \textbf{9},  
Walter de Gruyter \& Co., Berlin  (1988).


\bibitem{Io96a} D. Ioffe, {\em A note on the extremality of the disordered state for the Ising 
model on the Bethe lattice}, Lett. Math. Phys. \textbf{37}, 137-143 (1996).  


\bibitem{Io96b} D. Ioffe, {\em Extremality of the disordered state for the Ising 
model on  general trees}, Trees (Versailles), 3-14, Progr. Probab.  \textbf{40}, Birkh\"auser, Basel (1996). 

\bibitem{KeSt66}
H. Kesten,  , B. P.  Stigum, 
{\em Additional limit theorems for indecomposable multidimensional Galton-Watson processes}, 
Ann. Math. Statist. \textbf{37},  1463Ð1481 (1966).

\bibitem{Ly00}  R. Lyons, {\em 
Phase transitions on nonamenable graphs.  
In: Probabilistic techniques in equilibrium and nonequilibrium statistical physics, }  
J. Math. Phys. \textbf{41} no. 3, 1099--1126  (2000).


\bibitem{JaMo04}
S. Janson, E. Mossel, 
{\em Robust reconstruction on trees is determined by the second eigenvalue, } 
Ann. Probab. \textbf{32} no. 3B, 2630--2649 (2004).

\bibitem{Mo01} E. Mossel, {\em Reconstruction on trees: Beating the second eigenvalue}, 
Ann. Appl. Prob. \textbf{11} 285-300 (2001).  

\bibitem{MoPe03} E. Mossel, Y.Peres, {\em Information flow on trees}, 
The Annals of Applied Probability 
\textbf{13} no. 3, 817Ð844 (2003).

\bibitem{GeMo07} A. Gerschenfeld, A. Montanari,  
{\em Reconstruction for models on random graphs}, 
48th Annual IEEE Symposium on Foundations of Computer Science (FOCS), 194--204 (2007).

\bibitem{PePe06} R. Pemantle, Y. Peres, 
{\em The critical Ising model on trees, concave recursions and nonlinear capacity}, 
preprint, arXiv:math/0503137v2 [math.PR] (2006). 

\bibitem{Ma04} F. Martinelli, A. Sinclair, D. Weitz, {\em Fast mixing for independents sets, colorings and other model on trees}, 
Comm. Math. Phys. \textbf{250} no. 2, 301--334 (2004).

\bibitem{MeMo06} M. Mezard, A. Montanari, {\em Reconstruction on trees and spin glass transition}, 
 J. Stat. Phys. \textbf{124}  no. 6, 1317--1350 (2006).
 
 \bibitem{Bo06} C. Borgs, J.T. Chayes, E. Mossel, S.Roch, {\em The Kesten--Stigmun reconstruction bound is tight for roughly symmetric binary channels},
 47th Annual IEEE Symposium on Foundations of Computer Science (FOCS), 518-530 (2006).
 
 \bibitem{Sly08} A. Sly, {\em Reconstruction of symmetric Potts models},
 preprint, arXiv:0811.1208 [math.PR] (2008).
\end{thebibliography}
\end{document}